\documentclass[sn-mathphys-num]{sn-jnl}
\usepackage{amssymb}
\usepackage{amsmath}
\usepackage{tabularx}
\usepackage{bm}
\usepackage{booktabs}

\usepackage{derivative}

\usepackage{siunitx}

\DeclareSIUnit\cell{cell}
\DeclareSIUnit\day{day}

\renewcommand{\vec}[1]{\mathbf{#1}}

\begin{document}

\title[VI Solver for Simulating Burn Injury Healing]{Nondimensionalization, Variational Inequalities and Uncertainty Quantification in a Morphoelastic Model for Burn Injury Healing}

\author*[1,2]{\fnm{Wilbert} \spfx{den} \sur{Hertog}}\email{wilbert.denhertog@uhasselt.be}
\author[1,2]{\fnm{Fred} \sur{Vermolen}}\email{fred.vermolen@uhasselt.be}

\affil*[1]{\orgname{UHasselt – Hasselt University}, \orgdiv{Faculty of Sciences}, \orgaddress{\street{Agoralaan - Gebouw D}, \city{Diepenbeek}, \postcode{3590}, \country{Belgium}}}

\affil[2]{\orgname{UHasselt – Hasselt University}, \orgdiv{Data Science Institute (DSI)}, \orgaddress{\street{Agoralaan - Gebouw D}, \city{Diepenbeek}, \postcode{3590}, \country{Belgium}}}

\abstract{Mathematical modelling of burn injury healing involves coupling biological processes with mechanical deformation, resulting in computationally stiff, nonlinear partial differential equations that severely limit the scalability of existing numerical solutions. We present a computational framework addressing these challenges. By nondimensionalizing the model and performing scaling analysis, we show that chemotaxis and inertia are negligible, permitting a quasi-static momentum approach. To enforce non-negativity constraints, we implement a Variational Inequality formulation using the Finite Element Method in FEniCSx, and verify it via numerical convergence studies and mass balances. We identify high-order parameter interactions, specifically equilibrium collagen density and signalling dynamics, as the primary drivers of solver stability. A sensitivity analysis of nearly four million simulations reveals that isolated parameter variations predominantly govern the maximum wound contraction. In contrast, higher-order interactions determine the timing of maximum contraction and the permanent contracture. However, no significant pairwise parameter interactions were found, suggesting that the model's behaviour is governed by higher-order interactions, highlighting the inherent complexity of the wound-healing process.
}


\maketitle

\section{Introduction}
According to the WHO, burn injuries are a major public health problem, with an estimated \num{180000} deaths annually~\cite{who_burns_2023}, and many more injuries. Those who do not die from their (severe) injuries quite often suffer from long-term consequences, e.g. wound contraction and hypertrophic scar thickening. Mathematical modelling of burn injury healing can help understand the biological and mechanical processes that govern healing, can be used to develop treatments, and can be done without the need for costly (and unethical) in vivo experiments on animals.

The mathematical description of the evolution of skin after a (severe) burn injury, or wound healing in general, is a complex multiphysics challenge. It needs to connect microscopic biological processes, such as cell proliferation and differentiation, with macroscopic mechanical deformation of the skin. Different paradigms can be used to model tissue repair, e.g. agent-based models, cellular automata, or continuum models. The first two paradigms provide a microscopic description of the healing process; they are not ideal for modelling macroscopic skin deformation over larger portions of the human body. However, they can provide valuable insights into the fundamental processes of wound healing. Continuum models, on the other hand, are ideal for modelling macroscopic skin deformation and are thus well-suited for large-scale tissue modelling.

In this work, we present a new framework for numerical simulations of skin evolution following a severe burn injury. Our starting point is the existing mathematical model for burn injuries developed by Koppenol~\cite{koppenol_biomedical_2017} and Egberts~\cite{egberts_silico_2023}. This model is a coupled system of six nonlinear partial differential equations (PDEs) that models the concentrations and densities of fibroblasts, myofibroblasts, generic signalling molecules and collagen, as well as the deformation velocity and the strain of the skin. Biologically, these variables are coupled: fibroblasts proliferate and differentiate into myofibroblasts stimulated by signalling molecules. Myofibroblasts exert contractile forces that drive skin deformation, and they also deposit collagen, which affects the skin's mechanical properties. We use morphoelasticity theory to model the growth and deformation of the skin. The model has so far been explored through a custom Matlab implementation, which has provided an important foundation but is less suited to rapid scaling across multiple spatial dimensions and to large parameter-sensitivity analyses. In addition, the implementation has not yet been nondimensionalized, which makes it more difficult to assess the relative importance of the different parameters involved.

From a computational perspective, the existing implementation also presents several practical limitations in terms of extensibility and robustness. First, the reliance on dimensional parameters (that span more than 20 orders of magnitude) can lead to poorly conditioned matrices, which can cause numerical instabilities and convergence issues. Second, the use of fractional powers in the equations for fibroblasts, myofibroblasts, and the momentum equation requires that the biological variables remain strictly positive. Even a small numerical oscillation around zero could lead to negative values, rendering the fractional powers undefined and causing solver failure. Third, the custom code is less flexible when changing equations and boundary conditions, or when extending the model to multiple spatial dimensions. These considerations motivate a transition to a high-level automated finite element framework, such as FEniCSx.

The primary contribution of this work is the development of a highly efficient, scalable, and open-source computational framework that overcomes the mathematical and numerical limitations of previous implementations. Specifically, this framework introduces: (i) a nondimensionalization and scaling analysis of the existing mathematical model, which allows us to simplify the model by neglecting negligible terms, (ii) a reproducible open-source implementation of the simplified model in FEniCSx using a Variational Inequality formulation to enforce non-negativity constraints on the cell concentrations and chemical densities, (iii) a study of the solver's stability as a function of the parameters, and (iv) a large parameter sensitivity analysis using close to four million simulations. The sensitivity analysis reveals that isolated parameter variations influence the maximum wound contraction, whereas the day on which the maximum occurs and the final contraction are influenced by high-order parameter interactions.

The rest of this work is structured as follows. In Section~\ref{sect:model}, we present the mathematical model for the healing of burn injuries, and we perform a nondimensionalization and scaling analysis to simplify the model. In Section~\ref{sect:numerical_method}, we present the numerical method used to solve the equations. In Section~\ref{sect:results}, we present the results of the numerical simulations, including the convergence studies, the mass balance verification and the parameter sensitivity analyses. Finally, in Section~\ref{sect:conclusion}, we conclude this work and discuss future research directions.

\section{Mathematical Model}\label{sect:model}
In this section, we present the mathematical model for the healing of burn injuries as used by Egberts \cite{egberts_silico_2023} and Koppenol \cite{koppenol_biomedical_2017}. We then nondimensionalize the model and derive a simplified version based on a scaling analysis.
\subsection{Governing Equations}
The general, three-dimensional mathematical model for the healing of burn injuries developed and improved by Egberts \cite{egberts_silico_2023} and Koppenol \cite{koppenol_biomedical_2017} is given by the following equations:
\begin{align}
     & \mdv{z}{t}+z\nabla\cdot\vec{v}=-\nabla\cdot\vec{J}_z+R_z,\qquad z\in\{N,M,c,\rho\},\label{eq:chemicals}                                                                                                                         \\
     & \rho_t\left(\mdv{\vec{v}}{t}+\vec{v}(\nabla\cdot\vec{v})\right)=\nabla\cdot\bm{\sigma}+\vec{f},\label{eq:momentum}                                                                                                              \\
     & \bm\sigma = \mu_1\text{sym}(\nabla\vec{v})+\mu_2(\text{tr}(\text{sym}(\nabla\vec{v}))\bm{I})+\left(\frac{E\sqrt{\rho}}{1+\nu}\right)\left(\bm\varepsilon+\text{tr}(\bm\varepsilon)\left(\frac{\nu}{1-2\nu}\right)\bm{I}\right), \\
     & \mdv{\bm\varepsilon}{t}+\bm\varepsilon\text{skw}(\nabla\vec{v})-\text{skw}(\nabla\vec{v})\bm\varepsilon+(\text{tr}(\bm\varepsilon)-1)\text{sym}(\nabla\vec{v})=-\bm G,\label{eq:strain}
\end{align}
where $z=z(\vec{x},t)\in\mathbb{R},\ \vec{v} = \vec{v}(\vec{x},t)\in\mathbb{R}^3 \text{ and } \bm\varepsilon = \bm\varepsilon(\vec{x},t)\in\mathbb{R}^{3 \times 3}$. $\mdv{}{t}$ denotes the material derivative and $z\nabla\cdot\vec{v}$ represents passive convection. Here, $N$ and $M$ are the fibroblast and myofibroblast concentrations respectively (in \si{\cell\per\cm\cubed}), $c$ and $\rho$ are the signalling molecules and collagen densities (in \si{\gram\per\cm\cubed}). Further, $\vec{v}$ is the deformation velocity vector (in \si{\cm\per\day}), $\bm\sigma$ is the $3 \times 3$ stress tensor (in \si{\gram\per\cm\per\day\squared}), and $\bm\varepsilon$ is the nondimensional strain tensor. The system is closed by defining the fluxes $\vec{J}_z$, the reaction terms $R_z$ for each chemical and appropriate boundary and initial conditions. The transport of (myo)fibroblasts includes diffusion and chemotaxis towards the gradient of signalling molecules:
\begin{align}
    \vec{J}_N & = -D_F(N+M)\nabla N + \chi_FN\nabla c \\
    \vec{J}_M & = -D_F(N+M)\nabla M + \chi_FM\nabla c
\end{align}
The signalling molecules diffuse according to Fick's law, and for the collagen molecules, we assume no active transport:
\begin{align}
    \vec{J}_c    & = -D_c\nabla c \\
    \vec{J}_\rho & = \vec{0}
\end{align}
The reaction terms describe the proliferation, differentiation and decay of the chemicals:
\begin{align}
    R_N    & =r_F\left(1+\frac{r_F^\text{max}c}{a_c^I+c}\right)(1-\kappa_F(N+M))N^{1+q}-k_FcN-\delta_NN,                                      \\
    R_M    & =r_F\left(\frac{\left(1+r_F^\text{max}\right)c}{a_c^I+c}\right)(1-\kappa_F(N+M))M^{1+q}+k_FcN-\delta_MM,                         \\
    R_c    & =k_c\left(\frac{c}{a_c^{II}+c}\right)(N+\eta^IM)-\delta_c\frac{(N+\eta^{II}M)\rho}{1+a_c^{III}c}c,                               \\
    R_\rho & =k_\rho\left(1+\frac{k_\rho^{\text{max}}c}{a_c^{IV}+c}\right)(N+\eta^IM)-\delta_\rho\frac{(N+\eta^{II}M)\rho}{1+a_c^{III}c}\rho.
\end{align}
The body forces vector $\vec{f}$ and growth tensor $\bm G$ are given by:
\begin{equation}
    \vec{f}=\nabla\left(\frac{\xi M\rho}{R^2+\rho^2}\right), \quad \bm G = \zeta\frac{(N+\eta^{II}M)c}{1+a_c^{III}c}\bm\varepsilon.
\end{equation}

Typical dimensional parameter values and brief descriptions are provided in Appendix~\ref{app:dimensional_params}. For further information regarding the biology behind the model, the theses of Egberts \cite{egberts_silico_2023} and Koppenol \cite{koppenol_biomedical_2017} can be studied. An important note is that in the model used by Egberts and Koppenol, the unit of force is given in $\si{\newton}=\si{\kilo\gram\meter\per\second\squared}$. This unit is incompatible with the other units used in the model, namely \si{\gram}, \si{\centi\meter} and \si{\day}. Therefore, all parameters involving force have been converted to \si{\gram\centi\meter\per\day\squared} in this work. The conversion factor\footnote{$\qty{1}{\newton} = \qty{1}{\kilo\gram\meter\per\second\squared} = \qty{1000}{\gram}\cdot\qty{100}{\centi\meter}\cdot\left(\frac{1}{86400}\,\si{\day}\right)^{-2} = \qty{7.46496e14}{\gram\centi\meter\per\day\squared}$} is $\qty{1}{\newton}\approx\qty{7.46e14}{\gram\centi\meter\per\day\squared}$.

\subsection{Nondimensionalization and Scaling Analysis}
\label{subsect:nondim_scaling}
To improve the numerical implementation of the model and to study the relative importance of the different parameters, a nondimensionalization and scaling analysis is performed. We define the following characteristic scales:
\begin{gather*}
    x = L_c x^*, \quad t = t_c t^*, \\[1ex]
    N=N_cN^*, \quad M=N_cM^*, \quad c=c_cc^*, \quad \rho=\rho_c\rho^*, \\[1ex]
    \vec{v} = v_c \vec{v}^*, \quad \bm\sigma = \sigma_c\bm\sigma^*.
\end{gather*}
Note that the myofibroblast concentration $M$ is scaled with the characteristic fibroblast concentration $N_c$, since they share the same dimensions and order of magnitude. The reference scales used for nondimensionalization are shown in Table~\ref{table:scales}. Note that variables denoted with an asterisk ($^*$) are nondimensional. From these scaling parameters, we can derive the following nondimensional derivatives:
\begin{equation}
    \nabla(\cdot) = \frac{1}{L_c}\nabla^*(\cdot), \quad \nabla \cdot (\cdot) = \frac{1}{L_c}\nabla^* \cdot (\cdot), \quad \pdv{(\cdot)}{t} = \frac{1}{t_c}\pdv{(\cdot)}{t^*}, \quad \mdv{(\cdot)}{t} = \frac{1}{t_c}\mdv{(\cdot)}{t^*}.
\end{equation}
\begin{table}[htbp]
    \centering
    \caption{Reference scales used for nondimensionalization. The values are based on the typical parameter values shown in Table~\ref{table:parameters}, but should be recalculated if different parameter values are used.}
    \label{table:scales}
    \begin{tabular}{l l l}
        \toprule
        \textbf{Scale}                & \textbf{Definition}        & \textbf{Value/Unit}                         \\
        \midrule
        Length                        & $L_c = 1$                  & \qty{1}{\centi\meter}                       \\
        (Myo)fibroblast concentration & $N_c = \overline{N}$       & \qty{15000}{\cell\per\cm\cubed}             \\
        Collagen density              & $\rho_c = \overline{\rho}$ & \qty{0.125}{\gram\per\cm\cubed}             \\
        signalling molecules density  & $c_c = \tilde{c}$          & \qty{1e-8}{\gram\per\cm\cubed}              \\
        Time                          & $t_c = 1/\delta_N$         & \qty{52.6}{\day}                            \\
        Velocity                      & $v_c = L_c/t_c$            & \qty{0.019}{\cm\per\day}                    \\
        Stress                        & $\sigma_c = \mu_1 / t_c$   & \qty{1.42e15}{\gram\per\day\squared\per\cm} \\
        \bottomrule
    \end{tabular}
\end{table}

Substituting the scaling variables into the fibroblast equation \eqref{eq:chemicals} (here shown for $z=N$) and multiplying both sides of the equation by $t_c/N_c$ results in:
\begin{equation}
    \mdv{N^*}{t^*}+N^*(\nabla^*\cdot\vec{v}^*)=-\nabla^*\cdot\vec{J}^*_N+R_N^*,
\end{equation}
where
\begin{equation}
    \vec{J}^*_N = -D_F^*(N^*+M^*)\nabla^* N^* + \chi_F^* N^* \nabla^* c^*,
\end{equation}
and
\begin{equation}
    R_N^* = r_F^*\left(1+\frac{r_F^\text{max}c^*}{(a_c^I)^* + c^*}\right)(1-\kappa_F^* (N^* + M^*)) (N^*)^{1+q} - k_F^*c^*N^* - \delta_N^* N^*.
\end{equation}
Here, the following dimensionless parameters have been introduced:
\begin{gather*}
    D_F^*=\frac{D_FN_c t_c}{L_c^2}, \quad \chi_F^*=\frac{\chi_F c_c t_c}{L_c^2}, \quad r_F^* = r_F t_c N_c^q, \quad (a_c^I)^*=\frac{a_c^I}{c_c}, \\[1ex]
    \kappa_F^* = \kappa_F N_c, \quad k_F^* = k_F t_c c_c, \quad \delta_N^* = \delta_N t_c.
\end{gather*}
Applying an identical procedure to the myofibroblast, signalling molecules, and collagen equations results in the following nondimensional equations (where we have dropped the asterisks on the model parameters $N,M,c,\rho,\vec{v},x$ and $t$ and on the derivatives for readability):
\begin{align}
    \mdv{N}{t} + N (\nabla \cdot \vec{v})       & = -\nabla \cdot \left( -D_F^* (N+M) \nabla N + \chi_F^* N \nabla c \right) + R_N^*, \\
    \mdv{M}{t} + M (\nabla \cdot \vec{v})       & = -\nabla \cdot \left( -D_F^* (N+M) \nabla M + \chi_F^* M \nabla c \right) + R_M^*, \\
    \mdv{c}{t} + c (\nabla \cdot \vec{v})       & = -\nabla \cdot ( -D_c^* \nabla c ) + R_c^*, \label{eq:nondim_c}                    \\
    \mdv{\rho}{t} + \rho (\nabla \cdot \vec{v}) & = R_\rho^*,\label{eq:nondim_rho}
\end{align}
where
\begin{align}
    R_N^*    & = r_F^*\left(1+\frac{r_F^\text{max}c}{\left(a_c^I\right)^* + c}\right)(1-\kappa_F^* (N + M)) N^{1+q} - k_F^*cN - \delta_N^* N,               \\
    R_M^*    & = r_F^*\left(\frac{(1+r_F^\text{max})c}{\left(a_c^I\right)^*+c}\right)(1-\kappa_F^*(N+M))M^{1+q}+k_F^*cN-\delta_M^*M,                        \\
    R_c^*    & =k_c^*\left(\frac{c}{(a_c^{II})^*+c}\right)(N+\eta^IM)-\delta_c^*\frac{(N+\eta^{II}M)\rho}{1+(a_c^{III})^*c}c,                               \\
    R_\rho^* & =k_\rho^*\left(1+\frac{k_\rho^{\text{max}}c}{(a_c^{IV})^*+c}\right)(N+\eta^IM)-\delta_\rho^*\frac{(N+\eta^{II}M)\rho}{1+(a_c^{III})^*c}\rho.
\end{align}
All dimensionless parameters are shown in Table~\ref{table:dimensionless_params}. Performing the nondimensionalization of the momentum equation \eqref{eq:momentum} results in:
\begin{equation}
    \text{Re}\left(\mdv{\vec{v}}{t}+\vec{v}(\nabla\cdot\vec{v})\right)=\nabla\cdot\bm\sigma+\nabla\left(\frac{\xi^*M\rho}{(R^*)^2+\rho^2}\right),
\end{equation}
where
\begin{equation}
    \bm\sigma=\text{sym}(\nabla\vec{v})+\mu^*(\text{tr}(\text{sym}(\nabla\vec{v}))\bm{I})+\text{El}\frac{\sqrt{\rho}}{1+\nu}\left(\bm\varepsilon +\text{tr}(\varepsilon)\left(\frac{\nu}{1-2\nu}\right)\bm{I}\right).
\end{equation}
The nondimensional strain equation \eqref{eq:strain} becomes:
\begin{equation}
    \mdv{\bm\varepsilon}{t}+\bm\varepsilon\text{skw}(\nabla\vec{v})-\text{skw}(\nabla\vec{v})\bm\varepsilon+(\text{tr}(\bm\varepsilon)-1)\text{sym}(\nabla\vec{v}) = -\zeta^*\frac{(N+\eta^{II}M)c}{1+(a_c^{III})^*c}\bm\varepsilon.\label{eq:nondim_eps}
\end{equation}
\begin{table}[htbp]
    \centering
    \small
    \caption{Dimensionless parameters derived from the scaling analysis. Again, the values are based on the typical parameter values shown in Table~\ref{table:parameters}, but should be recalculated if different parameter values are used.}
    \label{table:dimensionless_params}
    \begin{tabular}{l l l | l l l}
        \toprule
        \textbf{Param}             & \textbf{Definition}        & \textbf{Value} & \textbf{Param} & \textbf{Definition}        & \textbf{Value} \\
        \midrule
        $\left(a_c^I\right)^*$     & $a_c^I / c_c$              & 1              & El             & $E\sqrt{\rho_c}/\sigma_c$  & 59.5           \\
        $\left(a_c^{II}\right)^*$  & $a_c^{II} / c_c$           & 1              & $k_c^*$        & $k_cN_ct_c/c_c$            & 23.7           \\
        $\left(a_c^{III}\right)^*$ & $a_c^{III}c_c$             & 2.05           & $k_F^*$        & $k_Ft_cc_c$                & 6              \\
        $\left(a_c^{IV}\right)^*$  & $a_c^{IV} / c_c$           & \num{0.1}      & $k_\rho^*$     & $k_\rho N_ct_c/\rho_c$     & \num{0.603}    \\
        $\chi_F^*$                 & $\chi_Fc_ct_c/L_c^2$       & \num{1.05e-9}  & $\kappa_F^*$   & $\kappa_FN_c$              & \num{0.015}    \\
        $D_c^*$                    & $D_ct_c/L_c^2$             & \num{0.171}    & $\mu^*$        & $\mu_2/\mu_1$              & 1              \\
        $D_F^*$                    & $D_FN_ct_c/L_c^2$          & \num{0.947}    & $r_F^*$        & $r_Ft_cN_c^q$              & 1.02           \\
        $\delta_c^*$               & $\delta_cN_c\rho_ct_c$     & \num{49.3}     & $R^*$          & $R/\rho_c$                 & 7.96           \\
        $\delta_M^*$               & $\delta_Mt_c$              & 3.16           & Re             & $\rho_tv_c^2/\sigma_c$     & \num{2.77e-19} \\
        $\delta_N^*$               & $\delta_Nt_c$              & 1              & $\xi^*$        & $\xi N_c/(\rho_c\sigma_c)$ & \num{2778.9}   \\
        $\delta_\rho^*$            & $\delta_\rho N_c\rho_ct_c$ & \num{0.603}    & $\zeta^*$      & $\zeta N_ct_cc_c$          & 3              \\
        \bottomrule
    \end{tabular}
\end{table}
In Table~\ref{table:dimensionless_params} we give the two new dimensionless numbers: the Reynolds number $\text{Re}=\frac{\rho_tv_c^2}{\sigma_c}=\frac{\rho_tv_cL_c}{\mu_1}$ and the elasticity number $\text{El}=\frac{E\sqrt{\rho_c}}{\sigma_c}$. The Reynolds number indicates the relative importance of inertial forces to viscous forces. The elasticity number represents the ratio between elastic forces and viscous forces.

There are four nondimensional numbers that stand out in Table~\ref{table:dimensionless_params}: $\chi_F^*=\num{1.05e-9}$, $\kappa_F^*=\num{0.015}$, $\text{Re}=\num{2.77e-19}$ and $\xi^*=\num{2778.9}$. The very low value of the nondimensional chemotactic parameter $\chi_F^*$, compared to the value of $D_F^*=\num{0.947}$, shows that chemotaxis has a negligible effect on steering the diffusion process of (myo)fibroblasts. The low value of the nondimensional crowding factor $\kappa_F^*$ indicates that changes in proliferation due to crowding are significant only at very high cell densities. Since the fibroblasts and myofibroblasts have been nondimensionalized, we expect that normally the term $\kappa_F^*(N+M)\ll 1$, but locally in the wounded area, this might not always be the case. Because of this, we choose to retain the crowding term. The very low value of the Reynolds number $\text{Re}=\num{2.77e-19}$ indicates that inertial forces are negligible compared to viscous forces, which justifies the use of a quasi-static approach for solving the momentum equation, which is a common assumption in biological applications. In the momentum equation, we encounter multiple scale differences. The nondimensional body force term $\xi^*=\num{2778.9}$ is much larger than the nondimensional elastic term $\text{El}=59.5$, which in turn is much larger than the nondimensional viscous term $\mu^*=1$. This means that the body forces are dominant over the elastic forces, which are dominant over the viscous forces. We chose to retain the viscous term in the model to preserve the skin's viscoelastic properties. However, we can expect that the viscous forces have a negligible effect on the solution compared to the body forces and elastic forces.

\subsection{Updated Mathematical Model}\label{subsect:updated_model}
Based on the scaling analysis in Section~\ref{subsect:nondim_scaling}, we can simplify the full mathematical model by neglecting the chemotactic terms in the (myo)fibroblast equations and by adopting a quasi-static approach for the momentum equation. We then get the following simplified equations for the (myo)fibroblast concentrations and the momentum equation:
\begin{align}
    \mdv{N}{t} + N (\nabla \cdot \vec{v})                                     & = \nabla \cdot \left(D_F^* (N+M) \nabla N \right) + R_N^*,\label{eq:nondim_N} \\
    \mdv{M}{t} + M (\nabla \cdot \vec{v})                                     & = \nabla \cdot \left(D_F^* (N+M) \nabla M \right) + R_M^*,\label{eq:nondim_M} \\
    \nabla\cdot\bm\sigma+\nabla\left(\frac{\xi^*M\rho}{(R^*)^2+\rho^2}\right) & =\vec{0},\label{eq:nondim_momentum}
\end{align}
where all variables and parameters are nondimensional. This model is closed by appropriate nondimensional boundary and initial conditions. The parameters are defined in Table~\ref{table:dimensionless_params}.

From here on, we will work with this simplified model, given by Equations~\eqref{eq:nondim_c}, \eqref{eq:nondim_rho}, \eqref{eq:nondim_eps} and \eqref{eq:nondim_N}--\eqref{eq:nondim_momentum}. We will also use the simplified (nondimensional) definitions of the flux and reaction terms, given by:
\begin{align}
    \vec{J}_N    & = D_F^*(N+M)\nabla N,\label{eq:flux_N} \\
    \vec{J}_M    & = D_F^*(N+M)\nabla M,\label{eq:flux_M} \\
    \vec{J}_c    & = D_c^*\nabla c,\label{eq:flux_c}      \\
    \vec{J}_\rho & = \vec{0},\label{eq:flux_rho}
\end{align}
and from now on, we drop the asterisks on the reaction terms $R_N^*$, $R_M^*$, $R_c^*$ and $R_\rho^*$ for readability.

\subsection{One-Dimensional Model Reduction}
The equations in Section~\ref{subsect:updated_model} are valid in three dimensions. However, to more easily facilitate parameter sensitivity analysis, we reduce the model to one spatial dimension. We take the spatial domain to be $x\in[0,L]$, where, in this paper, we set $L=\qty{10}{\centi\meter}$. Time is considered in the interval $t\in[0,T]$, where $T=\qty{365}{\day}$ in this study. In one dimension, the velocity vector $\vec{v}$ reduces to a scalar velocity $v$ and the strain tensor $\bm\varepsilon$ reduces to a scalar $\varepsilon$. The stress tensor $\bm\sigma$, the force vector $\vec{f}$ and the growth tensor $\bm G$ all reduce to (nondimensional) scalars $\sigma = \pdv{v}{x}+\text{El}\sqrt{\rho}\varepsilon$, $f=\pdv{}{x}\left(\frac{\xi^*M\rho}{(R^*)^2+\rho^2}\right)$ and $G=\zeta^*\frac{(N+\eta^{II}M)c}{1+\left(a_c^{III}c\right)}\varepsilon$. All gradients and divergences become derivatives in the $x$-direction. The one-dimensional version of the simplified model then becomes:
\begin{align}
    \pdv{N}{t} + \pdv{(Nv)}{x}                                           & = \pdv{}{x}\left(D_F^*(N+M)\pdv{N}{x}\right)+R_N,\label{eq:1d_N} \\
    \pdv{M}{t} + \pdv{(Mv)}{x}                                           & = \pdv{}{x}\left(D_F^*(N+M)\pdv{M}{x}\right)+R_M,\label{eq:1d_M} \\
    \pdv{c}{t} + \pdv{(cv)}{x}                                           & = D_c^* \pdv[2]{c}{x} + R_c,\label{eq:1d_c}                      \\
    \pdv{\rho}{t} + \pdv{(\rho v)}{x}                                    & = R_\rho,\label{eq:1d_rho}                                       \\
    \pdv{\sigma}{x} +f                                                   & =0,\label{eq:1d_force}                                           \\
    \pdv{\varepsilon}{t}+v\pdv{\varepsilon}{x}+(\varepsilon-1)\pdv{v}{x} & = -G.\label{eq:1d_strain}
\end{align}
For the force balance, Equation~\eqref{eq:1d_force}, we assumed $\nu=0$ and $\mu_2=0$, since there is no shear or bulk viscosity in one dimension.
To close the model, we need to define appropriate boundary and initial conditions. We assume that the left and right boundaries of the domain are both healthy skin, and that somewhere inside the domain (far away from the boundaries), there is a burn wound. This means that at the left and right boundaries, the (myo)fibroblast concentrations and the density of signalling molecules should take the values of healthy skin, given by $\overline{N}$, $\overline{M}$, and $\overline{c}$. Nondimensionally, this means: $N=1$, $M=0$ and $c=0$ at the boundaries. The deformation velocity $v$ is zero at the boundaries, since we assume that far away from the wound the displacement is zero. Furthermore, we do not have any boundary conditions for $\rho$ and $\varepsilon$ \cite{egberts_one-dimensional_2021}. For the initial conditions at time $t=0$, we assume that there is a burn injury in the centre of the domain, whose length is equal to $L_w=\qty{4}{\centi\meter}$. The initial wound length and the centre of the wound $x_c$ can be varied in the simulations. On the interface between the wound and healthy skin, we assume a smooth transition using a hyperbolic tangent indicator function $s(x)\in[0,1]$:
\begin{equation} s(x) = \frac{1}{2} \left[ 1 - \tanh\left( \alpha \left( |x - x_c| - \frac{L_w}{2} \right) \right) \right], \end{equation}
where $\alpha$ is a smoothing parameter (set to $\alpha=5$ in this study). The initial conditions for $0 < x < L$ are then given by:
\begin{align}
    N(x,0)           & = s(x)\tilde{N} + (1-s(x))\overline{N},       \\
    M(x,0)           & = 0                                           \\
    c(x,0)           & = s(x)\tilde{c}+(1-s(x))\overline{c},         \\
    \rho(x,0)        & = s(x)\tilde{\rho} + (1-s(x))\overline{\rho}, \\
    v(x,0)           & =0,                                           \\
    \varepsilon(x,0) & =0,
\end{align}
where $\tilde{N}$, $\tilde{c}$ and $\tilde{\rho}$ are the fibroblast concentration, signalling molecules density and collagen density, respectively in wounded skin and $\overline{N}$, $\overline{c}$ and $\overline{\rho}$ are the concentrations and densities in healthy skin. These values are given in Table~\ref{table:parameters} and should be nondimensionalized using the scales in Table~\ref{table:scales}.

In the simplified model, we neglect inertia (as explained in Section~\ref{subsect:updated_model}). To verify the validity of this neglect, we will compare the results of the full model (i.e., the model Egberts used) with those of the simplified model in Section~\ref{subsect:baseline}.

\section{Numerical Method}\label{sect:numerical_method}
Since an analytical solution to the model has not been found, we must solve it numerically. We solve this model monolithically using continuous Galerkin finite elements with linear Lagrangian basis functions. The time integration is performed using the backward Euler method, and the resulting nonlinear system is solved using a reduced-space active-set Newton method. The time step is adaptive and depends on the number of Newton iterations needed to solve the nonlinear system. The implementation is done in Python (version 3.13), using the DOLFINx library \cite{dolfinx2023}, together with UFL \cite{Alns2014} for the variational formulation and Basix \cite{Scroggs2022,Scroggs2022_1} for the mixed-element basis functions.

The weak and variational inequality formulations, as well as the time stepping strategy, are described in the following sections.
\subsection{Weak Formulation}
Koppenol presented the model's weak formulation in his PhD thesis \cite{koppenol_biomedical_2017}. We refer to his thesis for the details. In short, we multiply each equation by an appropriate test function, integrate over the domain and apply integration by parts where necessary. The superscripts $n$ and $n+1$ denote the previous and current time step, respectively (e.g.: $J_N^{n+1} = D_F^*(N^{n+1}+M^{n+1})\pdv{N^{n+1}}{x}$). This results in the weak formulation: Find $(N^{n+1},M^{n+1},c^{n+1},\rho^{n+1},v^{n+1},\varepsilon^{n+1})\in[H^1(0,L)]^6$ at a fixed time satisfying the (nondimensional) boundary conditions
\begin{gather*}
    N^{n+1}(0) = N^{n+1}(L) = 1, \\[1ex]
    M^{n+1}(0) = M^{n+1}(L) = c^{n+1}(0) = c^{n+1}(L) = v^{n+1}(0) = v^{n+1}(L) = 0.
\end{gather*}
such that for all test functions $(\phi_N,\phi_M,\phi_c,\phi_\rho,\phi_v,\phi_\varepsilon)\in V=[H^1_0(0,L)]^3 \times H^1(0,L) \times H^1_0(0,L) \times H^1(0,L)$:
\begin{align}
     & \int_0^L \left[\left(\frac{1}{\Delta t}(N^{n+1}-N^n) + \pdv{(N^{n+1}v^{n+1})}{x} - R_N^{n+1}\right)\phi_N + J_N^{n+1}\odv{\phi_N}{x}\right] \, dx = 0,\label{eq:wf_N}                                                    \\[1ex]
     & \int_0^L \left[\left(\frac{1}{\Delta t}(M^{n+1}-M^n) + \pdv{(M^{n+1}v^{n+1})}{x} - R_M^{n+1}\right)\phi_M + J_M^{n+1}\odv{\phi_M}{x}\right] \, dx = 0,\label{eq:wf_M}                                                    \\[1ex]
     & \int_0^L \left[\left(\frac{1}{\Delta t}(c^{n+1}-c^n) + \pdv{(c^{n+1}v^{n+1})}{x} - R_c^{n+1}\right)\phi_c + J_c^{n+1}\odv{\phi_c}{x}\right] \, dx = 0,\label{eq:wf_c}                                                    \\[1ex]
     & \int_0^L \left[\left(\frac{1}{\Delta t}(\rho^{n+1}-\rho^n) + \pdv{(\rho^{n+1} v^{n+1})}{x} - R_\rho^{n+1}\right)\phi_\rho\right] \, dx = 0,\label{eq:wf_rho}                                                             \\[1ex]
     & \int_0^L \left[\left(\sigma^{n+1}+\tilde{f}^{n+1}\right)\odv{\phi_v}{x}\right]\, dx = 0,\label{eq:wf_v}                                                                                                                  \\[1ex]
     & \int_0^L \left[\left(\frac{1}{\Delta t}(\varepsilon^{n+1}-\varepsilon^n) + v^{n+1}\pdv{\varepsilon^{n+1}}{x} + (\varepsilon^{n+1}-1)\pdv{v^{n+1}}{x} + G^{n+1}\right)\phi_\varepsilon\right] \, dx = 0\label{eq:wf_eps},
\end{align}
where $\tilde{f} = \frac{\xi^*M\rho}{(R^*)^2+\rho^2}$. The boundary terms from integration by parts vanish due to the chosen boundary conditions.
\subsection{Variational Inequality Framework}
The concentrations and densities of the (myo)fibroblasts, signalling molecules and collagen have to be non-negative. To enforce this constraint, we reformulate the problem as a Variational Inequality (VI) problem. We define the (closed and convex) feasible set of spatial solutions as
\begin{equation}
    \mathcal{K} = \{ \bm{\phi} \in [H^1(0,L)]^6 : \phi_i(x) \geq 0 \text{ for } i \in \{N,M,c,\rho\}, \;\forall x \in [0,L] \}.
\end{equation}
For any fixed time $t\in(0,T]$ we then seek a solution such that $(N(\cdot,t),M(\cdot,t),c(\cdot,t),\rho(\cdot,t),v(\cdot,t),\varepsilon(\cdot,t))\in\mathcal{K}$.
Rewriting the weak formulation from the previous subsection more compactly in its residual form results in: Find $\vec{u}^{n+1}=(N^{n+1}, M^{n+1}, c^{n+1}, \rho^{n+1}, v^{n+1}, \varepsilon^{n+1})$ such that
\begin{equation}
    F(\vec{u}^{n+1};\bm{\phi}) = 0 \quad \forall \bm{\phi} = (\phi_N, \phi_M, \phi_c, \phi_\rho, \phi_v, \phi_\varepsilon) \in V,\label{eq:weak_unconstrained}
\end{equation}
where $F(\vec{u};\bm{\phi}) = F_N(N,M,c,\rho,v,\varepsilon;\phi_N)+F_M(\cdot;\phi_M)+\cdots + F_\varepsilon(\cdot;\phi_\varepsilon)$ represents the sum of all residual forms from Equations \eqref{eq:wf_N}--\eqref{eq:wf_eps}. The constrained problem then becomes: Find $\vec{u}^{n+1}\in\mathcal{K}$ such that
\begin{equation}
    F(\vec{u}^{n+1};\vec{w}-\vec{u}^{n+1}) \geq 0 \quad \forall \vec{w}\in\mathcal{K}.\label{eq:variational_inequality}
\end{equation}
\subsubsection{Derivation of the Variational Inequality}
The VI formulation (Equation~\eqref{eq:variational_inequality}) arises from restricting the residual to only the feasible directions within $\mathcal{K}$. We derive the VI formulation using standard theory of variational inequalities (e.g. \cite{Kinderlehrer2000} and \cite{Glowinski_etc._1981}).

In the standard unconstrained setting, we seek solutions $\vec{u}_\text{unc}$ such that the residual $F(\vec{u}_\text{unc};\bm{\phi})$ is zero for all directions in the space $V$. But for the constrained solution $\vec{u}^*\in\mathcal{K}$, we can only test the feasible directions: directions that keep us inside $\mathcal{K}$. Otherwise, the solution could become negative. For any test function $\vec{w}\in\mathcal{K}$, the direction $\vec{w}-\vec{u}^*$ points in a direction that keeps us inside $\mathcal{K}$, since for any $\alpha\in[0,1]$, $\vec{u}^* + \alpha(\vec{w}-\vec{u}^*)$ is a convex combination of two points in $\mathcal{K}$ and is therefore also in $\mathcal{K}$. This leads to the VI formulation: we require that the residual $F(\vec{u}^*;\vec{w}-\vec{u}^*)$ is non-negative for all $\vec{w}\in\mathcal{K}$. We can verify this considering two cases:
\begin{itemize}
    \item \textbf{Case 1: Interior point.} Suppose that $\vec{u}^*$ is in the interior of $\mathcal{K}$ (i.e. all constrained variables are strictly positive), then we can move a small distance in any direction and still be in $\mathcal{K}$. This means that for a small enough $\epsilon > 0$ and any test function $\bm{\phi}\in V$, $\vec{u}^* \pm \epsilon\bm{\phi}$ is also in the interior of $\mathcal{K}$. Now, taking $\vec{w}=\vec{u}^* \pm \epsilon\bm{\phi}$ gives that $\vec{w}-\vec{u}^*=\pm\epsilon\bm{\phi}$. Substituting this into the VI formulation gives that $F(\vec{u}^*;\pm\epsilon\bm{\phi}) \geq 0$. Using the linearity of the residual form $F$ in its second argument we get that $\pm\epsilon F(\vec{u}^*;\bm{\phi}) \geq 0$, which forces $F(\vec{u}^*;\bm{\phi}) = 0$. So for interior points $\vec{u}^*$, this recovers the standard weak formulation.
    \item \textbf{Case 2: Boundary point.} As described above, at the boundary of $\mathcal{K}$ we can only test the directions that keep us inside $\mathcal{K}$. Specifically, whenever a constrained component of $\vec{u}^*$ is zero, we can only test directions that do not decrease that component (i.e., $\vec{w}-\vec{u}^*$ has to be nonnegative in that component). The condition $F(\vec{u}^*;\vec{w}-\vec{u}^*) \geq 0$ ensures that the residual `force' cancels out any variation into the forbidden region.
\end{itemize}
The VI formulation is equivalent to enforcing the following conditions for each constrained variable (here shown for $N$):
\begin{equation}
    N\geq0,\quad F_N(\vec{u}) \geq 0, \quad \text{and} \quad NF_N(\vec{u}) = 0 \text{ in } (0,L),
\end{equation}
where $F_N(\vec{u})$ is the nonlinear residual associated with the $N$-equation. These constraints are often referred to as the Karush-Kuhn-Tucker (KKT) conditions.

\subsubsection{Numerical Solution}
To solve the VI problem, we use PETSc's reduced-space active-set Newton method (\texttt{vinewtonrsls}) \cite{petsc-user-ref, Benson2006}. We access this solver through the \texttt{petsc4py} interface \cite{dalcinpazklercosimo2011}. This solver combines the following three strategies:
\begin{itemize}
    \item \textbf{Active sets}: At each Newton iteration, the solver checks which constraints are active (i.e. which constrained variables equal zero).
    \item \textbf{Reduced space Newton}: The Newton update is solved only for the inactive set (i.e. the unconstrained variables). The active set variables are kept fixed at zero.
    \item \textbf{Line search with bound enforcement}: After computing the Newton update, a line search is performed to find an appropriate step length. The step length is restricted such that the updated solution does not violate the constraints (i.e. does not become negative).
\end{itemize}
The algorithm terminates when the KKT conditions are satisfied to within a given tolerance.

\subsection{Adaptive Timestepping}\label{subsubsect:timestepping}
The timestep $\Delta t$ is adaptively adjusted based on the number of Newton iterations required to solve the nonlinear system at each time step. For the simulations done for this paper, we use the following strategy: we start with a small time step $\Delta t=\qty{0.05}{\day}$. We increase the time step by a factor of 1.2 if the number of Newton iterations needed to solve the system is less than 4, and we decrease the time step by a factor of 0.5 if the number of Newton iterations is greater than 8. If the number of iterations is between 4 and 8, we keep the time step unchanged. If, for a specific time step, the solver fails to converge within a maximum number of iterations (set to 50), we decrease the time step by multiplying it by 0.2 and try again. The minimum and maximum time step are set to $\qty{0.001}{\day}$ and $\qty{10}{\day}$ respectively. Note that the time step has to be scaled by the characteristic time $t_c$ to use in the solver. If the solver does not converge using the minimum time step, we terminate the simulation and report a failure to converge. This strategy allows us to use larger time steps when the solution is close to its equilibrium (e.g. after a few months of simulated wound healing), while using smaller time steps when the solution is changing rapidly (e.g. in the first few weeks of simulated wound healing).

\section{Results}\label{sect:results}
In this section, after showing the simulation results for the baseline parameters, we will first verify the implementation of the \texttt{vinewtonrsls} solver. We will do that by showing a time and space convergence analysis, and by verifying the mass balance of the (myo)fibroblasts, signalling molecules and collagen. Further, we compare the simplified model with the full model (i.e. the model Egberts and Koppenol used), to show that the simplifications we made are justified. Later, we will study the limitations of the solver, which results in a stability region for the parameters. Finally, we will perform a sensitivity analysis of the model's parameters to examine its biological outputs.

Visualization scripts for the figures in this section were generated with the assistance of large language models (Claude, Anthropic; Gemini, Google). The resulting visualizations were subsequently reviewed, validated and adapted by the authors. The data presented in these figures were generated by numerical simulations of the burn injury model, which were implemented and run by the authors.

\subsection{Baseline simulation}
\begin{figure}[htbp]
    \centering
    \includegraphics{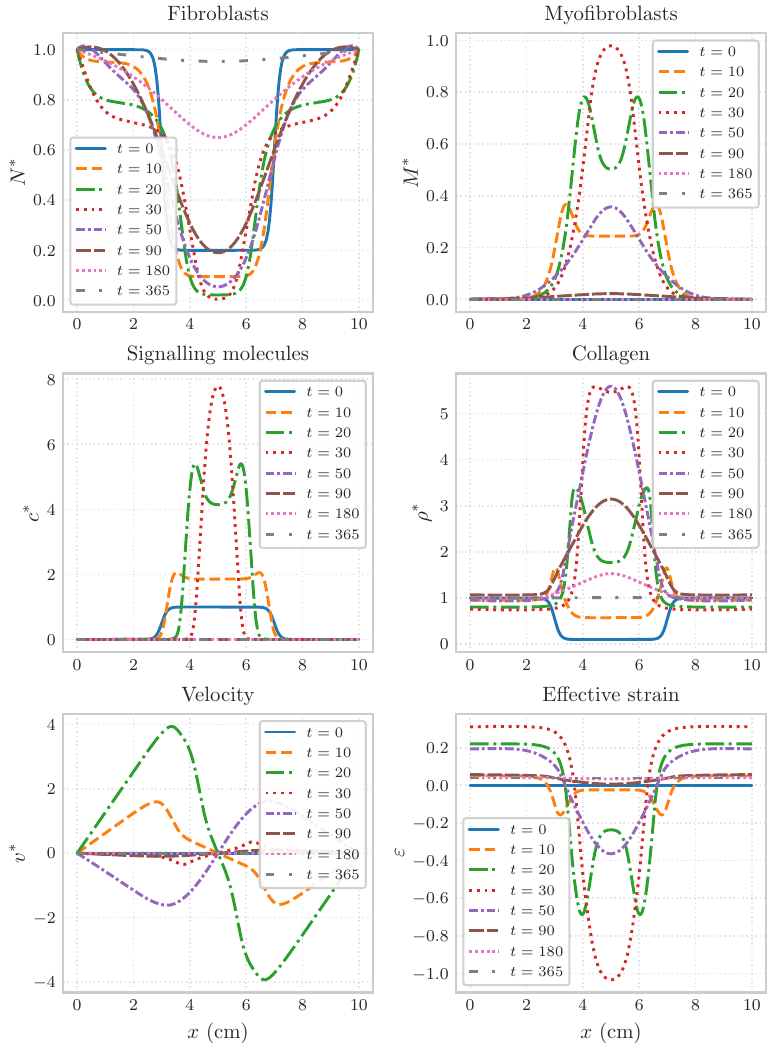}
    \caption{Spatial profiles of the nondimensional variables at different time points for the baseline simulation. The maximum contraction occurs around day 29. After that, the wound retracts.}
    \label{fig:baseline_simulation}
\end{figure}

\begin{figure}[htbp]
    \centering
    \includegraphics{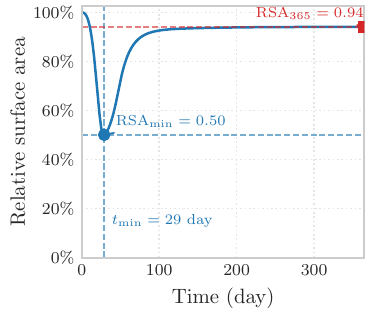}
    \caption{Relative surface area (RSA) of the wound over time for the baseline simulation. We see that the wound contracts rapidly in the first month, reaching an RSA of around 0.5 after 29 days. After that, the wound retracts relatively quickly, reaching an RSA of $0.9$ after 73 days, and then retracts more slowly, converging to a final RSA of around 0.94 (i.e., 6\% permanent contraction) after 1 year.}
    \label{fig:baseline_RSA}
\end{figure}

To get an idea of the model's output, we run a simulation using the baseline parameters from Table~\ref{table:parameters}. The spatial profiles of the nondimensional variables at different time points are shown in Figure~\ref{fig:baseline_simulation}. In this figure, we can see that the number of fibroblasts decreases (compared to the initial condition) across almost the entire domain during the first 90 days of injury healing. After that, it slowly converges to the equilibrium fibroblast concentration. The myofibroblasts and signalling molecules follow a similar pattern: during the first 30 days of wound healing, there is a large increase in myofibroblasts and signalling molecules, followed by a relatively quick convergence to the equilibrium solution, which is zero for both. The amount of collagen in the wound is initially low but increases locally to more than 5 times the equilibrium value during the first 30-50 days of wound healing. Afterwards, it converges to the equilibrium collagen density. The velocity and strain plots clearly show the contraction phase (until $t=\qty{30}{\day}$) and subsequent wound retraction. The relative surface area (RSA) of the wound is calculated by tracing the material points that are initially on the wound boundary (i.e. the points at $x=x_c \pm L_w/2$ at time $t=0$) and calculating the distance between those points at each time point, divided by the initial wound length $L_w$. The RSA over time is shown in Figure~\ref{fig:baseline_RSA}. The minimal RSA is approximately $0.5$, which means that, at the moment of maximum contraction, the wound length is approximately $50\%$ of its original length. That moment is after 29 days. This coincides with the moment the extreme values are reached for the (myo)fibroblasts, signalling molecules and collagen concentrations/densities and the strain, and also with the deformation velocity being close to zero over the whole domain at $t=\qty{30}{\day}$. The residual wound contraction after a year of simulated healing is approximately $6\%$.

\subsection{Convergence and Verification}
In this section, we verify the spatial and temporal convergence of the numerical solution, and then we will see that the solutions for the (myo)fibroblasts, signalling molecules and collagen also satisfy the mass balance. For all simulations done for the convergence and mass balance analysis, we use the following parameters: $L=\qty{10}{\centi\meter}$, $L_w=\qty{4}{\centi\meter}$, $x_c=\qty{5}{\centi\meter}$, and the parameter values from Table~\ref{table:parameters} (nondimensionalized using the scales from Table~\ref{table:scales}).

\subsubsection{Spatial and Temporal convergence}
Since we do not have an analytical solution to the model, we cannot compute the error of our numerical solution with respect to the true solution. Instead, we calculate the order of convergence by comparing the numerical solution with a reference solution computed on a very fine mesh and with a very small time step. The error is calculated using the $L^2$-norm, and the order of convergence is calculated by performing a least-squares fit on the $\log$-$\log$ plot of the error versus the mesh size $\Delta x$ (for spatial convergence) or time step $\Delta t$ (for temporal convergence). For the spatial convergence study, we use meshes with $N_x = \{10, 40, 160, 640, \num{2560}, \num{10240}\}$ elements and a fixed time step of $\Delta t = \qty{0.001}{\day}$. All simulations are run until $T=\qty{10}{\day}$. For the spatial reference solution, we use a mesh with $N_x = \num{40960}$ elements. For the temporal convergence study, we use a fixed mesh with $N_x=400$ elements and time steps $\Delta t \in [10^{-2}, 9.765625 \times 10^{-6}]$ day (6 levels, refined by factor 4). All simulations are run until $T=\qty{1}{\day}$. For the temporal reference solution, we use a time step of $\Delta t_{\text{ref}} = 2.44140625 \times 10^{-6}$ day. The timestepping for both convergence tests is done with a fixed time step, so the strategy described in Section~\ref{subsubsect:timestepping} is not used here.

\begin{figure}[htbp]
    \centering
    \includegraphics{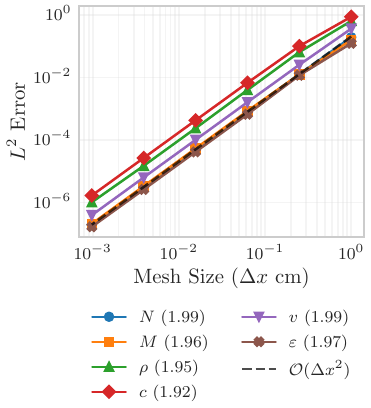}
    \caption{Spatial convergence plot. The orders of convergence are given in the legend between brackets.}
    \label{fig:spatial_convergence}
\end{figure}

\begin{figure}[htbp]
    \centering
    \includegraphics{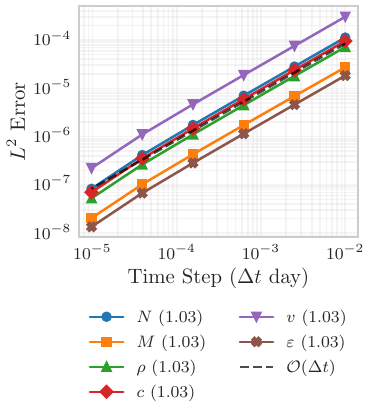}
    \caption{Temporal convergence plot. The orders of convergence are given in the legend between brackets.}
    \label{fig:temporal_convergence}
\end{figure}

Denoting the errors by $E_i$, and the corresponding mesh size by $\Delta x_i$ or time step by $\Delta t_i$, we fit a line to the points $(\log(\Delta x_i), \log(E_i))$ or $(\log(\Delta t_i), \log(E_i))$. The slope $p$ of this line is the order of convergence:
\begin{equation}
    \log(E_i) \approx p \log(h_i) + C,
\end{equation}
where $C$ is a constant and $h_i = \Delta x_i$ for spatial convergence and $h_i= \Delta t_i$ for the temporal convergence. The convergence plots are shown in Figures~\ref{fig:spatial_convergence} and \ref{fig:temporal_convergence}.

We find that the spatial convergence is second-order and the temporal convergence is first-order. This is consistent with the use of linear Lagrangian basis functions for spatial discretization and the backward Euler method for time integration. The entire convergence study (both spatial and temporal) takes approximately 3 hours on a single processor core.

\subsubsection{Mass balances}
\begin{figure}[htbp]
    \centering
    \includegraphics{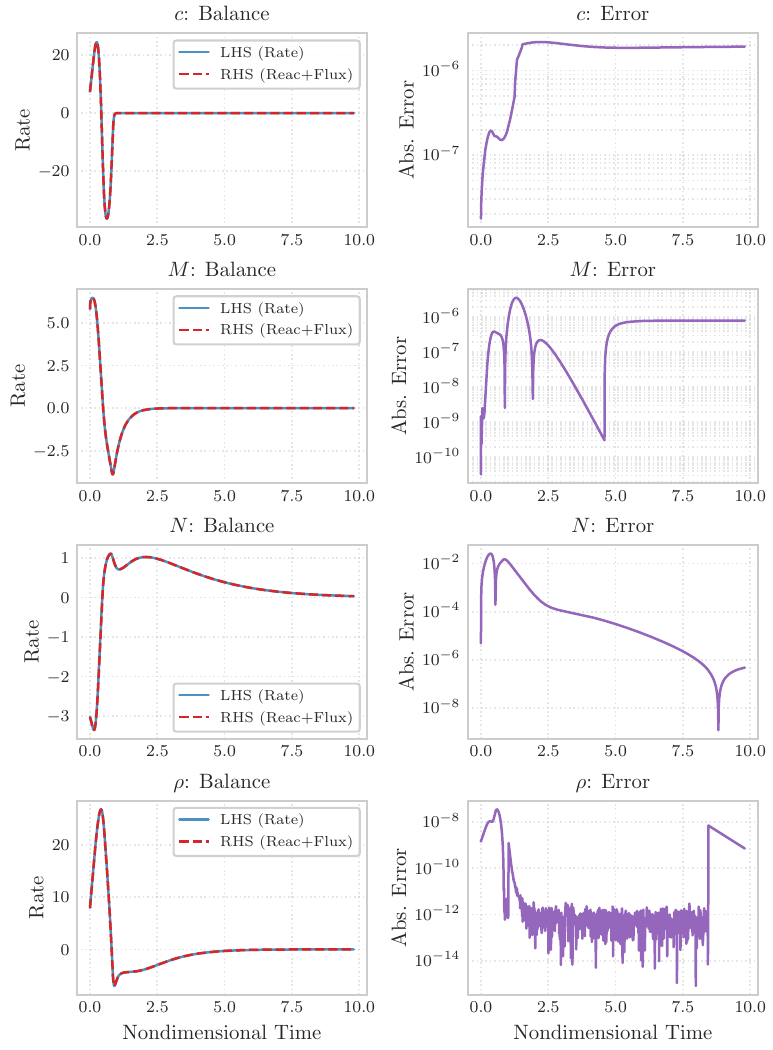}
    \caption{Mass balances. On the left, we show the change rate of the total mass of the (myo)fibroblasts, signalling molecules and collagen together with the sum of the reaction terms and the flux on the boundaries. On the right, we show the absolute mass-balance error.}
    \label{fig:mass_balance}
\end{figure}

The rate of change of the total mass of the (myo)fibroblasts, signalling molecules, and collagen should be equal to the sum of the reaction terms over the domain and the flux on the boundaries. Since $v=0$ on the boundary, there is no advective flux. So we only have diffusive fluxes on the boundaries. The mass balance is then given by:
\begin{equation}
    \frac{d}{dt}\int_0^L z\, dx = \int_0^L R_z\, dx + \left[ J_z \right]_0^L, \quad z \in \{N, M, c, \rho\} \label{eq:mass_balance}.
\end{equation}
In Figure~\ref{fig:mass_balance}, in the left column, we show the change rate of the total mass of the (myo)fibroblasts, signalling molecules and collagen, together with the sum of the reaction terms and the flux on the boundaries. In the right column, we show the absolute error of the mass balance, which is given by the absolute difference between the left and right-hand side of the mass balance equation (Equation~\eqref{eq:mass_balance}). For $M$, $c$, and $\rho$, we observe that the mass balance is satisfied up to an error of $\mathcal{O}(10^{-6})$, which is very good. For $N$, we observe a larger error of $\mathcal{O}(10^{-2})$. A possible explanation is that the boundary flux for $N$ is $\mathcal{O}(1)$, the same order as its reaction term, which results in spatial discretization errors in the boundary flux contributing significantly to the mass balance error. The other variables have maximal boundary fluxes that are $\mathcal{O}(10^{-3})$ or smaller, while their reaction terms are $\mathcal{O}(1)$ to $\mathcal{O}(10)$. We found that the mass balance error for $N$ is insensitive to the time step size but decreases proportionally with mesh refinement, consistent with the explanation that the error is dominated by the spatial discretization of the boundary flux. Overall, these results show that the mass balance is satisfied with very good accuracy, verifying the correctness of the solver implementation.

\subsection{Solver Stability Analysis}\label{subsect:solver_stability}
To study the possibilities of the solver, we verify its stability (i.e. whether the solver converges or not) by performing a sensitivity analysis. Apart from general execution errors, there are two main ways in which the solver (gracefully) fails: either the solver needs more than some maximum number of total Newton iterations to converge (where we sum the number of Newton iterations needed at each time step), or the solver needs a time step smaller than the minimum time step (as described in Section~\ref{subsubsect:timestepping}).

We perform the sensitivity analysis using the Sobol method~\cite{Sobol2001}. This is a variance-based method that decomposes the output variance (i.e., whether the solver fails) into contributions from each input parameter and their interactions. The total Sobol index $S_T$ of a parameter represents the total contribution of that parameter to the variance of the output, including both its direct effect (which is called the first-order effect $S_1$) and its interactions with other parameters (higher-order effects). We combined the Sobol method with Saltelli's sampling scheme~\cite{Saltelli2002}. One of the drawbacks of the Sobol method is that it needs a large number of model evaluations to accurately compute Sobol indices, namely $N_s\cdot(2(p+1))$ model evaluations, where $N_s$ is the number of samples and $p$ is the number of parameters. Note that $N_s$ should be a power of 2 for the sampling scheme to work. In our case, we varied $p=29$ parameters: the wound length and the parameters from Table~\ref{table:parameters} (except for $\overline{c}$, $\chi_F$, $k_\rho$, $\mu_2$, $\overline{M}$, $\nu$, $q$, $\rho_t$). We used $N_s=\num{16384}$ samples, which results in a total of \num{983040} model evaluations. Each model evaluation corresponds to a single simulation of the burn injury model, in which healing is simulated over one year. The parameters were sampled over the range of 0.5 to 1.5 times their typical values (as given in Table~\ref{table:parameters}).

\begin{figure}[htbp]
    \centering
    \includegraphics{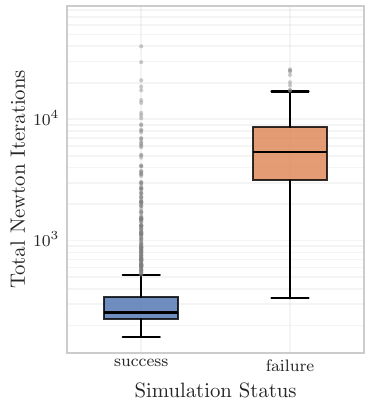}
    \caption{Distribution of Total Newton Iterations for successful and failed simulations. Successful simulations have a median of 255 and a mean of 667 iterations, while failed simulations have a median of \num{5355} and a mean of \num{6657} iterations. This indicates that the distributions are well separated, despite some overlap in the tails. Based on these results, we chose a cutoff of \num{6000} Newton iterations for the sensitivity analysis, retaining 98.4\% of the successful simulations and terminating 43.9\% of the runs that would eventually fail.}
    \label{fig:newton_iterations}
\end{figure}

Since each model evaluation can take from less than a second to several minutes, we first performed a preliminary analysis to determine a cutoff for the total number of Newton iterations over all time steps of a single simulation, which allows us to terminate simulations that are likely to fail early on, to save computational resources. The results of this preliminary analysis are shown in Figure~\ref{fig:newton_iterations}. Based on this analysis, we chose a cutoff of \num{6000} total Newton iterations (using a minimal time step of $\qty{0.001}{\day}$), which retains 98.4\% of the successful simulations and terminates 43.9\% of the runs that would eventually fail. We found that whenever a simulation reaches \num{6000} Newton iterations, it has an 83\% probability of ultimately failing. After applying this cutoff, we ran the sensitivity analysis, which took approximately 11.5 hours on 128 processing cores (1 load balancing core and 127 calculation cores) on the VSC HPC cluster. On average, each simulation took approximately 5.3 seconds.

\subsubsection{Sobol analysis}
\begin{figure}[htbp]
    \centering
    \includegraphics{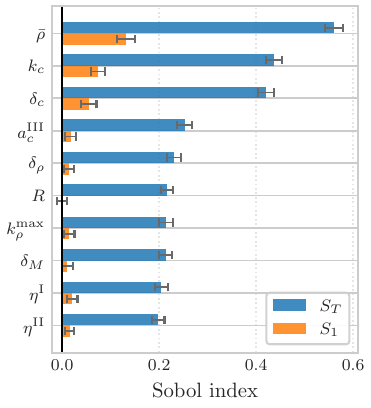}
    \caption{The ten most impactful parameters that influence the stability of the solver. $S_T$ is the total Sobol index, which represents the total contribution of a parameter to the variance of the output (i.e. whether the solver fails or not). $S_1$ is the first-order Sobol index, which represents the contribution of a parameter to the variance of the output on its own (i.e. without interactions with other parameters). The error bars represent the 95\% confidence intervals of the Sobol indices.}
    \label{fig:sobol_analysis}
\end{figure}

The Sobol analysis was performed on the model's binary (success/failure) output. The results of the Sobol analysis are shown in Figure~\ref{fig:sobol_analysis}. The ten most impactful parameters influencing the solver's stability are shown, sorted by their total Sobol index $S_T$. The three most impactful parameters are the equilibrium collagen density $\overline{\rho}$, the signalling molecule secretion rate $k_c$ and the signalling molecule degradation rate $\delta_c$, with $S_T$ scores of 0.56, 0.44 and 0.42, respectively. Notably, for each parameter, its first-order Sobol index $S_1$ is very low relative to its total Sobol index, indicating that parameter interactions are very important for the solver's stability. We have also examined the pairwise interactions between parameters (i.e., the second-order Sobol indices), but we have not found any significant interactions. This suggests that the interactions between parameters that influence the solver's stability are of higher order (i.e., involving three or more parameters).

\subsubsection{Dominant parameters}
\begin{figure}[htbp]
    \centering
    \includegraphics{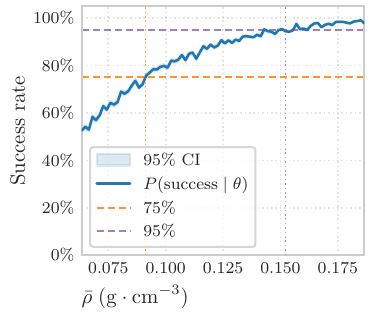}
    \caption{Convergence rate of the solver as a function of the equilibrium collagen density $\overline{\rho}$. The blue line shows the probability of a successful simulation for specific values of $\overline{\rho}$. The plot shows that the solver converges more slowly for low values of $\overline{\rho}$.}
    \label{fig:rho_eq}
\end{figure}

\begin{figure}[htbp]
    \centering
    \includegraphics{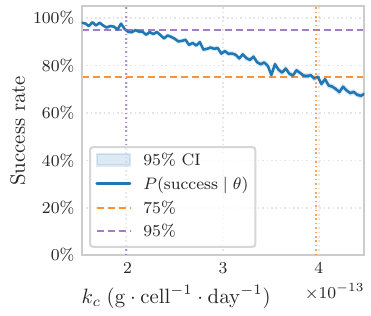}
    \caption{Convergence rate of the solver as a function of the signalling molecules secretion rate $k_c$. The blue line shows the probability of a successful simulation for specific values of $k_c$. The plot shows that the solver converges more slowly for higher values of $k_c$.}
    \label{fig:k_c}
\end{figure}

\begin{figure}[htbp]
    \centering
    \includegraphics{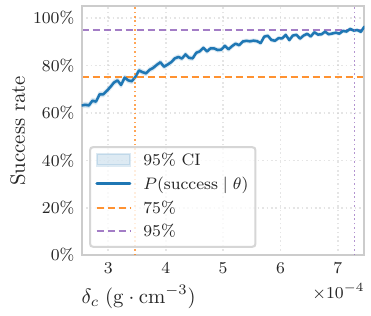}
    \caption{Convergence rate of the solver as a function of the signalling molecules breakdown rate $\delta_c$. The blue line shows the probability of a successful simulation for specific values of $\delta_c$. The plot shows that the solver converges more slowly for low values of $\delta_c$.}
    \label{fig:delta_c}
\end{figure}

For the three most impactful parameters, additional analysis was conducted to examine how the solver's convergence rate changes as a function of these parameters. For each of these three parameters, we calculated the percentage of successful simulations for specific parameter values across all simulations in the Sobol analysis. The results are shown in Figures~\ref{fig:rho_eq}, \ref{fig:k_c} and \ref{fig:delta_c}. These plots show that it is harder for the solver to converge for low values of $\overline{\rho}$, high values of $k_c$ and low values of $\delta_c$. This suggests there are two physical regimes that negatively impact the solver's stability: one in which the equilibrium collagen density is low, and one in which the amount of signalling molecules is high (since high $k_c$ and low $\delta_c$ both lead to more signalling molecules).

\subsubsection{Analysis of nondimensional parameters}
Since this model has dimensional constraints on $q$ and $k_\rho$ (which must be derived from other dimensional parameters), it is not feasible to perform a direct Sobol analysis on the nondimensional parameters. Furthermore, the Saltelli sampling scheme requires independently sampled inputs; however, the nondimensional parameters are nonlinear combinations of the same dimensional parameters and are therefore not independent. An example: $\overline{\rho}=\rho_c$, the characteristic collagen concentration, which is used in the calculation of six nondimensional parameters simultaneously: $\delta_c^*$, $\delta_\rho^*$, El, $k_\rho^*$, $R^*$ and $\xi^*$. This suggests that the high Sobol index of $\overline{\rho}$ may, at least in part, reflect its role as a scaling parameter. To complement the Sobol analysis of the dimensional parameters, we compute the nondimensional parameters for each of the \num{983040} simulations and calculate the Spearman correlation coefficient between each nondimensional parameter and the binary success/failure output. An important note in this analysis is that $\delta_N^*=\delta_Nt_c=\frac{\delta_N}{\delta_N}=1$ by construction, so it is not included in the correlation analysis. Further omitted are $\overline{N}$, $\overline{\rho}$ and $\tilde{c}$, since nondimensionally they also always equal one.

\begin{figure}[htbp]
    \centering
    \includegraphics{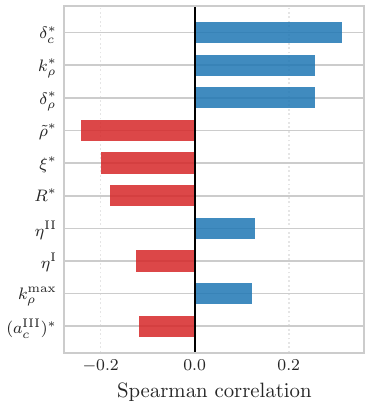}
    \caption{Spearman correlation coefficient between the nondimensional parameters and the binary success/failure output. The blue bars indicate a positive correlation; the red bars, a negative correlation. A positive (negative) correlation means that higher values of the parameter are associated with a higher (lower) probability of successful convergence. The parameters are sorted by the absolute value of their correlation coefficient.}
    \label{fig:correlation_nondim}
\end{figure}

The results of this analysis are shown in Figure~\ref{fig:correlation_nondim}. The highest-ranking parameters are $\delta_c^*$, $k_\rho^*$, and $\delta_\rho^*$, which are related to signalling molecules and collagen dynamics, consistent with the Sobol analysis results. A major difference in the results compared to the Sobol analysis is that the nondimensional initial collagen concentration $\tilde{\rho}^*$ in the wounded area is ranked fourth, whereas it was ranked 27th (out of 29 parameters) in the Sobol analysis. We also computed the distance correlation~\cite{Szkely2007} between the nondimensional parameters and the binary success/failure output and found rankings identical to those of the Spearman correlation, suggesting that the relationships between the nondimensional parameters and the solver's stability are mostly monotonic.

\subsection{Comparison of the baseline run}\label{subsect:baseline}
\begin{figure}[htbp]
    \centering
    \includegraphics{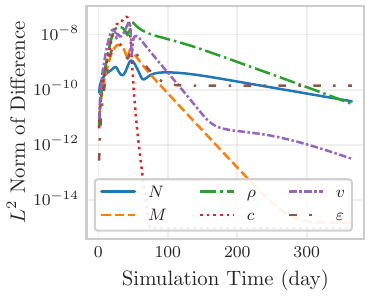}
    \caption{Results of the model comparison. The plot shows the $L^2$ norm of the nondimensional difference between the models per variable.}
    \label{fig:simulation_comparison_L2}
\end{figure}

To validate the new computational framework, we compare the simplified model's results with those of the original model (i.e., the one-dimensional version used by Egberts and Koppenol, which includes chemotactic and inertia terms). We run both models with the same baseline parameters from Table~\ref{table:parameters} and compare the results in this section.
For this comparison, we chose a fixed timestep of $\qty{0.01}{\day}$. At every timestep, we calculated the $L^2$ norm of the difference between the nondimensional variables in those two models. Since all variables are nondimensional and $\mathcal{O}(1)$, the $L^2$ norm of the difference directly represents relative errors. In Figure~\ref{fig:simulation_comparison_L2}, we see a maximum error of $\mathcal{O}(10^{-7})$, which is consistent with the prediction in Section~\ref{subsect:nondim_scaling} that $\chi_F^*$ and Re would have very little impact. Simulations with parameters other than the baseline also show very small errors, further supporting this claim.

\subsection{Parameter Sensitivity Analysis}\label{subsect:sensitivity}
Egberts studied the sensitivity of the model to its parameters by performing a one-at-a-time (OAT) sensitivity analysis, in which each parameter is varied one at a time while keeping the other parameters fixed at their baseline values~\cite{Egberts2021}. These results, however, cannot be directly compared with ours, since a different scaling was used in Egberts' sensitivity analysis (parameter variations of $\pm0.25$\% instead of $\pm25$\%). Since we now believe Egberts' scaling is unrealistic, we present the results of our sensitivity analysis without comparison to Egberts' results.

We performed a global sensitivity analysis using the Sobol method. This allows us to study both the first-order effects of each parameter (i.e., the effect of a parameter on the output without considering interactions with other parameters) and the total effects (i.e., the effect of a parameter on the output including interactions with other parameters). The baseline parameters from Table~\ref{table:parameters} are used in our study. The total parameter space is defined by varying each parameter over the range 0.75 to 1.25 times its baseline value. The sampling is performed using Saltelli's scheme, which enables efficient computation of the Sobol indices. For this sensitivity analysis, we focus on clinically relevant output variables that can be compared to medical data. Specifically, we study the following model outputs: the minimum relative surface area (RSA) reached during the simulation, the time at which it is reached, and the final RSA after 1 year of simulated healing. We found that for the Sobol analysis, we need $N_s=\num{65536}$ to obtain acceptable confidence intervals, resulting in a total of $\num{65536}\cdot (2\cdot 29 + 2) = \num{3932160}$ model evaluations. These simulations were run on 128 processing cores (1 load-balancing core and 127 calculation cores) on the VSC HPC cluster and took approximately 17 hours. So, on average, each simulation took approximately 2 seconds. It is important to note that \num{57858} simulations (1.47\%) failed to converge, likely because the parameters were in an unstable region. To calculate the Sobol indices, we need numerical model outputs for each simulation. Because of the block structure of the sampling scheme, we cannot simply discard the failed simulations; we would need to discard entire blocks of simulations, which would lead to an approximate 59\% data loss (with a failure rate of 1.47\%, the chance that a block of 60 simulations contains at least one failed simulation is approximately 59\%). Because of this, we decided to impute the failed simulation results with the median value of the successful simulations for each output variable. This does influence the results of the sensitivity analysis; however, we expect the most important parameters to still be identified as such, since the percentage of failed simulations is relatively low. In Section~\ref{subsubsect:min_rsa} we present the results of the minimal RSA and the time at which that minimum is reached, and in Section~\ref{subsubsect:rsa_end} we present the results of the permanent contraction (i.e. the final RSA after one year of simulated healing). The table of Sobol indices for all parameters is in Appendix~\ref{app:sensitivity_scores}.

\subsubsection{Maximal contraction}\label{subsubsect:min_rsa}
\begin{figure}[htbp]
    \centering
    \includegraphics{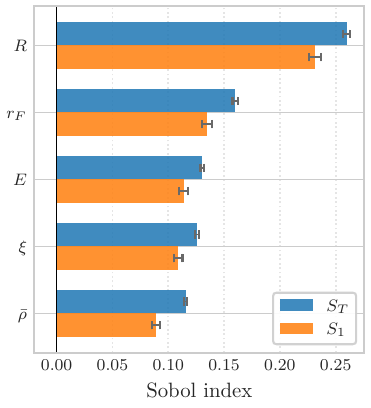}
    \caption{top five Sobol indices for the minimal RSA, sorted by total Sobol index $S_T$. The error bars represent the 95\% confidence intervals of the Sobol indices.}
    \label{fig:min_rsa}
\end{figure}

\begin{figure}[htbp]
    \centering
    \includegraphics{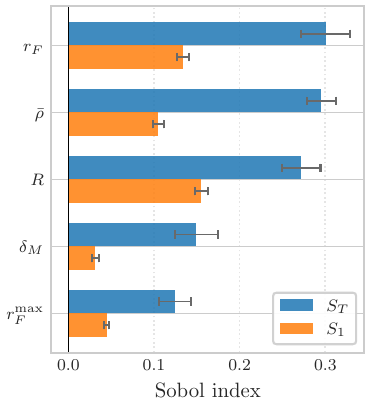}
    \caption{top five Sobol indices for the day at which the minimal RSA is reached, sorted by total Sobol index $S_T$. The error bars represent the 95\% confidence intervals of the Sobol indices.}
    \label{fig:min_rsa_day}
\end{figure}

In Figures~\ref{fig:min_rsa} and \ref{fig:min_rsa_day} we show the top five Sobol indices for the minimal RSA and the day at which that minimum is reached, sorted by total Sobol index $S_T$. Noteworthy is that the $S_1$ and $S_T$ values are very close for the minimal RSA (Figure \ref{fig:min_rsa}), which suggests that the maximum contraction is mostly controlled by the direct effects of the parameters. In contrast, for the day of minimal RSA (Figure \ref{fig:min_rsa_day}), the $S_1$ values generally are not even 50\% of the $S_T$ values, which suggests that more than half of the influence of the parameters on the day of minimal RSA is due to interactions between parameters. There is some overlap in the top five parameters that influence the minimal RSA and the day of minimal RSA, namely, $R$, $r_F$, and $\overline{\rho}$ appear in both plots. The differences are that the minimal RSA is also strongly influenced by $E$ and $\xi$, while the day of minimal RSA is, in turn, strongly influenced by $\delta_M$ and $r_F^\text{max}$.

\subsubsection{Permanent contraction}\label{subsubsect:rsa_end}
\begin{figure}[htbp]
    \centering
    \includegraphics{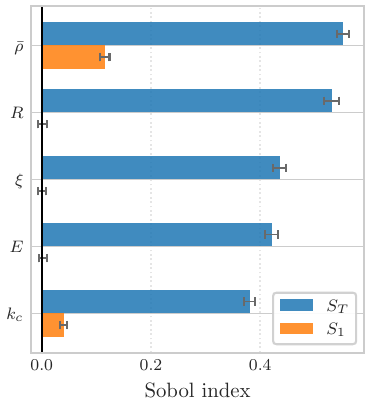}
    \caption{top five Sobol indices for the final RSA after one year of simulated healing, sorted by total Sobol index $S_T$. The error bars represent the 95\% confidence intervals of the Sobol indices.}
    \label{fig:sobol_rsa_end}
\end{figure}

The top five Sobol indices for the final RSA are shown in Figure~\ref{fig:sobol_rsa_end}. The most influential parameters are $\overline{\rho}$, $R$, $\xi$, $E$ and $k_c$. We see that $S_1$ is much smaller than $S_T$ for all parameters, which suggests that the interactions between parameters are very important for the permanent contraction and that no single parameter has a dominant influence on the permanent contraction on its own.

\section{Conclusion and discussion}\label{sect:conclusion}
This work presents a new computational framework for the mathematical modelling of burn injury healing. We build on the mathematical model of Koppenol~\cite{koppenol_biomedical_2017} and Egberts~\cite{egberts_silico_2023}. We nondimensionalized the model, enabling a scaling analysis that showed the chemotactic and inertia terms are negligible compared to the other terms. The numerical implementation of the model is done in FEniCSx, using the Variational Inequality framework to enforce non-negativity constraints. We verified the implementation of the solver by performing a convergence analysis and by verifying the mass balances of the (myo)fibroblasts, signalling molecules and collagen. We also performed two large-scale analyses: a solver stability analysis using nearly 1 million simulations, which revealed that solver stability is governed by high-order parameter interactions, and a global sensitivity analysis of the biological model outputs using nearly 4 million simulations, which showed that permanent wound contraction is driven by interactions between both biological and mechanical parameters. We discuss the main findings and their implications in the following paragraphs.

\subsection{Nondimensionalization and scaling analysis}
To our knowledge, this is the first time that a nondimensionalization and scaling analysis of the burn injury model has been presented. Koppenol briefly mentions the characteristic scales he used for the nondimensionalization, but does not present the nondimensionalized equations. The scaling analysis showed that two terms can be safely neglected. The Reynolds number $\text{Re} = \num{2.77e-19}$ shows that inertia can be neglected; we can use a quasi-static force balance instead of a momentum equation. The chemotactic coefficient $\chi_F^* = \num{1.05e-9}$ is roughly nine orders of magnitude smaller than the diffusion coefficient $D_F^* = \num{0.947}$, which shows that chemotaxis can be neglected compared to diffusion. These numbers are specific to the baseline parameters, but their orders of magnitude will be similar for other parameter values. We verified both simplifications numerically: the comparison of the baseline simulation with the original model (which includes the chemotactic and inertia terms) showed a maximum error of $\mathcal{O}(10^{-7})$.

However, in reality, chemotaxis towards increasing concentrations of signalling molecules plays a key role in the healing of (burn) injuries~\cite{Singer1999}, so in future work it would be interesting to study how chemotaxis could be incorporated into the model. The nondimensional viscosity ($\mu^*=1$) was deliberately retained, despite being small compared to the nondimensional body force term ($\xi^*=\num{2778.9}$), to preserve the viscoelastic nature of the problem.

\subsection{Variational Inequality solver}
The numerical implementation of the model is done in FEniCSx, using the Variational Inequality framework to enforce non-negativity constraints for the chemicals. In the previous MATLAB implementation, small numerical oscillations could lead to negative concentrations. This causes issues with the fractional-power terms $N^{1+q}$ and $M^{1+q}$ in the fibroblast and myofibroblast equations, respectively, as well as with $\sqrt{\rho}$ in the force balance equation. $q$ is a negative number (for the baseline run, $q=-0.43$), which causes another numerical subtlety, namely that $1+q < 1$, and therefore the derivative of $N^{1+q}$ and $M^{1+q}$ with respect to $N$ and $M$ respectively, becomes unbounded as $N$ and $M$ approach zero. Since the solver uses a Newton method, this will lead to numerical issues when $N$ and $M$ approach zero. To avoid this, the implementation uses the regularization $u^\alpha\approx(u+\epsilon)^\alpha - \epsilon^\alpha$. Here, $u\in\{N,M,\rho\}$, and we use a small $\epsilon > 0$. This regularization preserves the exact zero at $u=0$, recovers $u^\alpha$ for $u\gg\epsilon$, and yields a bounded derivative for $u\geq 0$.

Besides the above-mentioned implementation details, the use of the FEniCSx framework also allows for easily modifying weak forms, boundary conditions, and initial conditions without the need to write assembly routines from scratch. FEniCSx also natively supports spatial parallelization, which is not yet used in the one-dimensional case, but will be very useful for future extensions of the model to higher dimensions.

Verification of the implementation was performed through spatial and temporal convergence studies, confirming second-order spatial convergence and first-order temporal convergence, consistent with the use of linear Lagrangian basis functions for spatial discretization and the backward Euler method for time integration. We also verified the mass balances of the (myo)fibroblasts, signalling molecules and collagen, which are satisfied up to an error of $\mathcal{O}(10^{-6})$ for $M$, $c$ and $\rho$, and up to an error of $\mathcal{O}(10^{-2})$ for $N$. The larger error for $N$ is due to the fact that the boundary flux of $N$ is $\mathcal{O}(1)$, the same order as its reaction term, which causes the spatial discretization errors in the boundary flux to contribute significantly to the mass balance error. We found that this mass balance error decreases proportionally with mesh refinement.

\subsection{Solver stability}
The solver stability study was done using nearly a million simulations. The results showed that the solver stability is governed almost entirely by high-order parameter interactions, since for every parameter, the first-order Sobol index $S_1$ is very low compared to the total Sobol index $S_T$. Also, no significant pairwise interactions have been found. This means that the solver stability cannot be reliably predicted or avoided by modifying individual parameters. The three most impactful parameters that influence the stability of the solver are the equilibrium collagen density $\overline{\rho}$, the signalling molecule secretion rate $k_c$ and the signalling molecules breakdown rate $\delta_c$. The solver struggles with low values of $\overline{\rho}$, and high concentrations of signalling molecules (high $k_c$ and low $\delta_c$). A possible explanation is that high concentrations of signalling molecules lead to large reaction terms in the equations, resulting in rapidly changing solutions that are harder to numerically solve. Low equilibrium collagen density $\overline{\rho}$ leads to another problem: decoupling of the strain and stress tensors. This coupling is present in the definition of $\sigma = \pdv{v}{x}+\text{El}\sqrt{\rho}\varepsilon$. If $\rho\approx0$, then $\sigma$ would just contain spatial variations in velocity, without any feedback from the strain tensor $\varepsilon$. The closer $\overline{\rho}$ gets to zero, the closer the force balance equation (Equation~\eqref{eq:1d_force}) will be to this regime, which then leads to numerical issues.

The Spearman correlation analysis also shows that the nondimensional parameters related to the signalling molecules and collagen dynamics ($\delta_c^*$, $k_\rho^*$ and $\delta_\rho^*$) are the most influential parameters for the stability of the solver. A notable difference from the dimensional Sobol analysis is that the nondimensional initial collagen concentration $\tilde{\rho}^*$ in the wounded area is ranked fourth, whereas it was ranked 27th (out of 29 parameters) in the dimensional analysis. This is likely due to the differences between the Sobol and Spearman analyses. The Sobol analysis captures the variance in model output, which is extremely low for $\tilde{\rho}$, with respect to the solver stability, while the Spearman analysis just captures the monotonic relationship between $\tilde{\rho}$ and the solver stability, which is high. The reason why the nondimensional initial collagen concentration $\tilde{\rho}^*$ is negatively correlated with the stability of the solver is due to the fact that $\tilde{\rho}^*=\frac{\tilde{\rho}}{\overline{\rho}}$, so higher values of $\tilde{\rho}^*$ correspond to lower values of $\overline{\rho}$, which we already know lead to a less stable solver.

The results obtained in this work can be connected to the linear stability conditions derived by Egberts~\cite{Egberts2021stability}:  $\delta_c\overline{\rho}\geq\frac{k_c}{a_c^{II}}$ and $q\delta_N\leq\kappa_Fr_F\overline{N}^{1+q}$. Note that the second constraint always holds for all simulations done in both the solver stability analysis and the parameter sensitivity analysis, since $q$ is negative in the whole parameter space, making the left-hand side of the second constraint strictly negative. For the dataset generated by the solver stability analysis (from Section~\ref{subsect:solver_stability}), we found that the first constraint is violated in $15.35\%$ of the simulations, and the crash rate conditional on the violation of this constraint is $48.57\%$ compared to $9.24\%$ when the constraint is satisfied. The dataset generated by the biological parameter sensitivity analysis (from Section~\ref{subsect:sensitivity}) shows a similar pattern. The first constraint is violated in $0.63\%$ of the simulations. The crash rate conditional on the violation of this constraint is $25.94\%$ compared to $1.32\%$ when the constraint is satisfied. This suggests that the linear stability conditions can be a useful practical tool for identifying problematic parameter regimes, but they are not sufficient to guarantee stability.

\subsection{Biological sensitivity analysis}
The global sensitivity analysis for the biological model outputs, based on nearly four million simulations, shows the following main results. The minimal RSA is mostly determined by the direct effects of the parameters, whereas the day of minimal RSA and the final RSA are mostly determined by parameter interactions. This is what we would expect, since the minimal RSA is primarily determined by mechanics: e.g., if the skin is stiff (high stiffness coefficient $E$), there simply cannot be as much contraction. Consequently, a lower skin stiffness would directly lead to more contraction. In contrast, the timing of peak contraction is more interaction-dominated, suggesting that it is inherently sensitive to the interaction between biological and mechanical processes. This is an expected result, as the speed of wound healing in patients is influenced by many factors, such as age and overall health. These factors are reflected in the model parameters through changing the rates of the biological processes (e.g., the proliferation rate of (myo)fibroblasts, the secretion and breakdown rates of signalling molecules, etc.).

For patients, the final RSA is the most relevant output, since it represents the permanent consequences of the burn injury. The finding that $S_1 \ll S_T$ for all parameters means that no single parameter has a dominant influence on permanent contraction on its own. The top five most influential parameters ($\overline{\rho}$, $R$, $\xi$, $E$ and $k_c$) span both mechanical and biological processes. This suggests that treatments targeting only one of these processes will not be the most effective, and that a combination targeting both will likely be needed to reduce permanent contraction.

Notably, none of the model outputs showed significant pairwise interactions between parameters, indicating that the interactions driving the day of minimal RSA and the final RSA are of higher order (i.e., involving three or more parameters), further highlighting the complexity of the healing process.

\subsection{Limitations and future work}
The current work has limitations that could be addressed in future work.

The time integration is based on the backward Euler method, which is a first-order method. For one-dimensional simulations, this is sufficient, as we can simply increase the temporal resolution to obtain accurate results without incurring excessive additional compute time. However, for extensions of the model to higher dimensions, higher-order time integration methods would be a natural improvement. A requirement for such a method is that it retains compatibility with adaptive timestepping, since in the first months of simulated healing, the solutions change rapidly, but after about half a year, they stabilize, allowing a larger time step.

The parameter ranges in the sensitivity analysis for the biological outputs were set to 0.75-1.25 times the baseline values, a somewhat arbitrary choice. Some parameters might have a much wider or narrower range of medically relevant values, or should not be varied at all. An improvement for future work would be to determine the medically relevant parameter ranges based on literature and expert knowledge, and to use those ranges in the sensitivity analysis.

Relatedly, the model parameters of the baseline are based on medical papers, as introduced by Koppenol~\cite{koppenol_biomedical_2017}. However, the model has not been fitted to patient data. So, for now, it is still impossible to say how well the model can predict the healing of real burn injuries. Model validation would be required before it can be used for clinical predictions or treatment planning.

The current model uses infinitesimal strain theory. For the large deformations that occur in severe burn injuries, a finite strain formulation would be more appropriate. This could be addressed in future work.

Finally, extending the framework to higher dimensions is now relatively straightforward. The one-dimensional reduction was adopted here to facilitate large-scale parameter studies, but the ultimate goal is to develop a model for patient-specific predictions, which would require a two- or three-dimensional model. For modelling wound depth (and scar thickening), we cannot use a full Eulerian description of the problem, since the domain itself changes over time. Instead, we would need to use a formulation that also tracks the domain, such as an Arbitrary Lagrangian-Eulerian (ALE) or phase-field formulation. The extension to higher dimensions would also allow for more complex geometries.

\section*{Statements and Declarations}

\subsection*{Acknowledgements}
The resources and services used in this work were provided by the VSC (Flemish Supercomputer Center), funded by the Research Foundation -- Flanders (FWO) and the Flemish Government.

\subsection*{Competing interests}
The authors declare no competing interests.

\subsection*{Code and data availability}
The FEniCSx implementation of the 1D model and the simulation datasets generated during this study are available in Zenodo at \url{https://doi.org/10.5281/zenodo.20595706}.

\subsection*{Funding}
This research received no specific grant from any funding agency in the public, commercial, or not-for-profit sectors.

\subsection*{Ethics approval}
Not applicable.

\appendix
\clearpage
\section{Dimensional Parameters}\label{app:dimensional_params}
\begin{table}[htbp]
    \centering
    \caption{Typical parameter values and their units used in the burn injury model.}
    \label{table:parameters}
    \begin{tabularx}{\textwidth}{ l r l X }
        \toprule
        \textbf{Parameter}      & {\textbf{Value}}     & \textbf{Unit}                                                & \textbf{Description}                                                   \\
        \midrule
        $a_{c}^{I}$             & $1\mathrm{e}{-8}$    & \si{\gram\per\cubic\centi\meter}                             & (Myo)fibroblast proliferation half-saturation constant                 \\
        $a_{c}^{II}$            & $1\mathrm{e}{-8}$    & \si{\gram\per\cubic\centi\meter}                             & Signalling molecules secretion half-saturation constant                \\
        $a_{c}^{III}$           & $2.05\mathrm{e}{8}$  & \si{\cubic\centi\meter\per\gram}                             & Degradation inhibition constant                                        \\
        $a_{c}^{IV}$            & $1\mathrm{e}{-9}$    & \si{\gram\per\cubic\centi\meter}                             & Collagen secretion half-saturation constant                            \\
        $\overline{c}$          & $0$                  & \si{\gram\per\cubic\centi\meter}                             & Signalling molecules in healthy skin                                   \\
        $\tilde{c}$             & $1\mathrm{e}{-8}$    & \si{\gram\per\cubic\centi\meter}                             & Signalling molecules in wounded skin                                   \\
        $\chi_F$                & $2\mathrm{e}{-3}$    & \si{\centi\meter\tothe{5}\per\gram\per\day}                  & Chemotactic coefficient                                                \\
        $D_c$                   & $3.25\mathrm{e}{-3}$ & \si{\centi\meter\squared\per\day}                            & Signalling molecules diffusion coefficient                             \\
        $D_F$                   & $1.2\mathrm{e}{-6}$  & \si{\centi\meter\tothe{5}\per\cell\per\day}                  & (Myo)fibroblast diffusion coefficient                                  \\
        $\delta_c$              & $5\mathrm{e}{-4}$    & \si{\centi\meter\tothe{6}\per\cell\per\gram\per\day}         & Breakdown rate signalling molecules                                    \\
        $\delta_M$              & $6\mathrm{e}{-2}$    & \si{\per\day}                                                & Breakdown rate myofibro-blasts                                         \\
        $\delta_N$              & $1.9\mathrm{e}{-2}$  & \si{\per\day}                                                & Breakdown rate fibroblasts                                             \\
        $\delta_{\rho}$         & $6.11\mathrm{e}{-6}$ & \si{\centi\meter\tothe{6}\per\cell\per\gram\per\day}         & Breakdown rate of collagen                                             \\
        $\eta^I$                & $2.0$                & {--}                                                         & Relative signalling secretion rate                                     \\
        $\eta^{II}$             & $0.45$               & {--}                                                         & Relative MMP secretion rate                                            \\
        $E$                     & $2.39\mathrm{e}{17}$ & \si{\gram\tothe{0.5}\centi\meter\tothe{0.5}\per\day\squared} & Stiffness coefficient ($E\sqrt{\rho} = \text{Young's Modulus}$)        \\
        $k_c$                   & $3\mathrm{e}{-13}$   & \si{\gram\per\cell\per\day}                                  & Signalling molecules secretion rate                                    \\
        $k_F$                   & $1.14\mathrm{e}{7}$  & \si{\centi\meter\cubed\per\gram\per\day}                     & Cell differentiation rate                                              \\
        $k_\rho$                & $9.54\mathrm{e}{-8}$ & \si{\gram\per\cell\per\day}                                  & Collagen secretion rate ($k_\rho=\delta_\rho(\overline{\rho})^2$)      \\
        $k_{\rho}^{\text{max}}$ & $10.0$               & {--}                                                         & Maximum collagen secretion rate enhancing factor                       \\
        $\kappa_F$              & $1\mathrm{e}{-6}$    & \si{\cubic\centi\meter\per\cell}                             & Crowding factor term                                                   \\
        $\mu_1$                 & $7.46\mathrm{e}{16}$ & \si{\gram\per\day\per\centi\meter}                           & Shear viscosity                                                        \\
        $\mu_2$                 & $7.46\mathrm{e}{16}$ & \si{\gram\per\day\per\centi\meter}                           & Bulk viscosity                                                         \\
        $\overline{M}$          & $0$                  & \si{\cell\per\cubic\centi\meter}                             & Myofibroblast concentration in healthy skin                            \\
        $\nu$                   & $0.49$               & {--}                                                         & Poisson's ratio                                                        \\
        $\overline{N}$          & $1.5\mathrm{e}{4}$   & \si{\cell\per\cubic\centi\meter}                             & Fibroblast concentration in healthy skin                               \\
        $\tilde{N}$             & $3\mathrm{e}{3}$     & \si{\cell\per\cubic\centi\meter}                             & Fibroblast concentration in wounded skin ($0.2\overline{N}$)           \\
        $q$                     & $-0.43$              & {--}                                                         & $q=(\log\delta_N-\log(r_F(1-\kappa_F \overline{N})))/\log\overline{N}$ \\
        $r_F$                   & $1.222$              & \si{\centi\meter\tothe{3q}\cell\tothe{-q}\per\day}           & Cell division rate                                                     \\
        $r_{F}^{\text{max}}$    & $2.0$                & {--}                                                         & Maximum cell division enhancing factor                                 \\
        $\overline{\rho}$       & $0.125$              & \si{\gram\per\cubic\centi\meter}                             & Collagen density in healthy skin                                       \\
        $\tilde{\rho}$          & $0.0125$             & \si{\gram\per\cubic\centi\meter}                             & Collagen density in wounded skin ($0.1\overline{\rho}$)                \\
        $\rho_t$                & $1.09$               & \si{\gram\per\cubic\centi\meter}                             & Skin density                                                           \\
        $R$                     & $0.995$              & \si{\gram\per\centi\meter\cubed}                             & Fixed constant                                                         \\
        $\xi$                   & $3.28\mathrm{e}{13}$ & \si{\gram\squared\per\day\squared\per\cell\per\centi\meter}  & Body force term                                                        \\
        $\zeta$                 & $3.8\mathrm{e}{2}$   & \si{\centi\meter\tothe{6}\per\cell\per\gram\per\day}         & Rate of morphoelastic change                                           \\
        \bottomrule
    \end{tabularx}
\end{table}
\newpage
\section{Sensitivity scores of all parameters}\label{app:sensitivity_scores}

\begin{table}[htbp]
    \centering
    \caption{
        Sobol sensitivity indices $S_1$ and $S_T$ for all parameters and three model outputs. Parameters are sorted by $S_T$ of the final RSA in descending order.
    }
    \label{table:sobol_sensitivity}
    \begin{tabular}{l r r r r r r}

        \toprule

                              & \multicolumn{2}{c}{\textbf{Min RSA}} & \multicolumn{2}{c}{\textbf{Time of}} & \multicolumn{2}{c}{\textbf{Final RSA}}                          \\
                              & \multicolumn{2}{c}{}                 & \multicolumn{2}{c}{\textbf{min RSA}} & \multicolumn{2}{c}{}                                            \\

        \cmidrule(lr){2-3}
        \cmidrule(lr){4-5}
        \cmidrule(lr){6-7}

        \textbf{Parameter}
                              & $S_1$                                & $S_T$
                              & $S_1$                                & $S_T$
                              & $S_1$                                & $S_T$                                                                                                  \\

        \midrule
        $\overline{\rho}$     & 0.089                                & 0.116                                & 0.105                                  & 0.295 & 0.115  & 0.551 \\
        $R$                   & 0.232                                & 0.260                                & 0.155                                  & 0.272 & 0.001  & 0.531 \\
        $\xi$                 & 0.109                                & 0.126                                & 0.039                                  & 0.089 & 0.001  & 0.436 \\
        $E$                   & 0.114                                & 0.130                                & 0.031                                  & 0.070 & 0.002  & 0.421 \\
        $k_c$                 & 0.039                                & 0.060                                & 0.001                                  & 0.106 & 0.040  & 0.381 \\
        $r_F$                 & 0.135                                & 0.160                                & 0.134                                  & 0.301 & -0.010 & 0.358 \\
        $\delta_c$            & 0.022                                & 0.042                                & 0.005                                  & 0.086 & 0.033  & 0.331 \\
        $\delta_\rho$         & 0.017                                & 0.027                                & 0.030                                  & 0.113 & 0.018  & 0.295 \\
        $\overline{N}$        & 0.001                                & 0.008                                & 0.011                                  & 0.055 & 0.006  & 0.275 \\
        $k_\rho^{\text{max}}$ & 0.015                                & 0.023                                & 0.027                                  & 0.111 & 0.015  & 0.271 \\
        $\delta_M$            & 0.069                                & 0.086                                & 0.031                                  & 0.150 & -0.003 & 0.265 \\
        $r_F^{\text{max}}$    & 0.047                                & 0.060                                & 0.045                                  & 0.125 & -0.005 & 0.260 \\
        $L_w$                 & 0.001                                & 0.010                                & 0.012                                  & 0.054 & -0.004 & 0.240 \\
        $a_c^{II}$            & 0.005                                & 0.011                                & 0.000                                  & 0.040 & -0.001 & 0.215 \\
        $a_c^{III}$           & 0.007                                & 0.020                                & 0.003                                  & 0.072 & 0.003  & 0.214 \\
        $\tilde{c}$           & 0.004                                & 0.011                                & 0.001                                  & 0.015 & 0.001  & 0.200 \\
        $\tilde{N}$           & 0.006                                & 0.009                                & 0.024                                  & 0.053 & 0.000  & 0.197 \\
        $\eta^I$              & 0.005                                & 0.014                                & 0.017                                  & 0.096 & 0.003  & 0.195 \\
        $k_F$                 & 0.000                                & 0.006                                & 0.002                                  & 0.033 & -0.003 & 0.191 \\
        $a_c^I$               & 0.007                                & 0.009                                & 0.007                                  & 0.031 & -0.003 & 0.169 \\
        $\delta_N$            & 0.001                                & 0.004                                & 0.004                                  & 0.034 & -0.004 & 0.168 \\
        $\mu_1$               & 0.000                                & 0.003                                & 0.000                                  & 0.008 & -0.001 & 0.165 \\
        $\eta^{II}$           & 0.003                                & 0.012                                & 0.003                                  & 0.063 & 0.003  & 0.123 \\
        $D_F$                 & 0.000                                & 0.002                                & 0.001                                  & 0.013 & -0.002 & 0.114 \\
        $\zeta$               & 0.001                                & 0.003                                & 0.001                                  & 0.019 & 0.013  & 0.105 \\
        $\tilde{\rho}$        & 0.000                                & 0.001                                & 0.000                                  & 0.009 & -0.002 & 0.071 \\
        $D_c$                 & 0.000                                & 0.001                                & -0.001                                 & 0.018 & -0.001 & 0.046 \\
        $a_c^{IV}$            & 0.000                                & 0.001                                & 0.000                                  & 0.014 & -0.002 & 0.037 \\
        $\kappa_F$            & 0.000                                & 0.001                                & 0.000                                  & 0.007 & 0.000  & 0.034 \\
        \bottomrule
    \end{tabular}
\end{table}
\clearpage
\bibliography{references}
\end{document}